\documentclass[11pt]{article}
\usepackage{psfig}
\def\X{{\cal X}}
\def\Y{{\cal Y}}
\def\C{{\cal C}}

\def\beq{\begin{eqnarray}}
\def\eeq{\end{eqnarray}}

\begin{document}

\fontsize{11}{14.5pt}\selectfont

\vspace*{0.6in}

\begin{center}

{\small Technical Report No.\ 0304,
 Department of Statistics, University of Toronto}

\vspace*{0.9in}

{\LARGE Markov Chain Sampling for Non-linear State Space Models\\[5pt]
        Using Embedded Hidden Markov Models }
\\[16pt]

{\large Radford M. Neal}\\[3pt]
 Department of Statistics and Department of Computer Science \\
 University of Toronto, Toronto, Ontario, Canada \\
 \texttt{http://www.cs.utoronto.ca/$\sim$radford/} \\
 \texttt{radford@stat.utoronto.ca}\\[10pt]

 30 April 2003
\end{center}

\vspace{8pt} 

\noindent \textbf{Abstract.}  I describe a new Markov chain method for
sampling from the distribution of the state sequences in a non-linear
state space model, given the observation sequence.  This method
updates all states in the sequence simultaneously using an embedded
Hidden Markov model (HMM).  An update begins with the creation of a
``pool'' of $K$ states at each time, by applying some Markov chain
update to the current state.  These pools define an embedded HMM whose
states are indexes within this pool.  Using the forward-backward
dynamic programming algorithm, we can then efficiently choose a state
sequence at random with the appropriate probabilities from the
exponentially large number of state sequences that pass through states
in these pools.  I show empirically that when states at nearby times
are strongly dependent, embedded HMM sampling can perform better than
Metropolis methods that update one state at a time.

\section{\hspace*{-7pt}Introduction}\label{sec-intro}\vspace*{-10pt}

Consider a state space model with observations $y_0,\ldots,y_{n-1}$, each
in some set $\Y$, and hidden states $x_0,\ldots,x_{n-1}$, each in some
set $\X$.  Suppose we know the dynamics of hidden states and the
observation process for this model.  Our task is to sample from the
distribution for the hidden state sequence given the observations.  

If the state space, $\X$, is finite, of size $K$, so that this is a
Hidden Markov Model (HMM), a hidden state sequence can be sampled by a
well-known forward-backwards dynamic programming procedure in time
proportional to $nK^2$.  Scott (2002) reviews this algorithm and
related methods.  If $\X = \Re^p$ and the dynamics and observation
process are linear, with Gaussian noise, an analogous adaptation of
the Kalman filter can be used.  For more general models, one might use
Markov chain sampling.  For instance, one could perform Gibbs sampling
or Metropolis updates for each $x_t$ in turn.  Such simple Markov
chain updates may be very slow to converge, however, if the states at
nearby times are highly dependent.

In this note, I describe how Markov chain sampling for these models
can be facilitated by using updates that are based on temporarily
embedding an HMM whose finite state space is a subset of $\X$, and
then applying the efficient HMM sampling procedure.

\section{\hspace*{-7pt}The Embedded HMM Algorithm}\label{sec-alg}\vspace*{-10pt}

In describing the algorithm, model probabilities will be denoted by
$P$ (which will denote probabilities or probability densities without
distinction, as appropriate for the state space, $\X$, and observation
space, $\Y$).  The initial state distribution is given by $P(x_0)$,
transition probabilities are given by $P(x_t\,|\,x_{t-1})$, and
observation probabilities are given by $P(y_t\,|\,x_t)$.  Our goal is
to sample from the conditional distribution $P(x_0,\ldots,x_{n-1}\,|\,
y_0,\ldots,y_{n-1})$, which we will abbreviate to
$\pi(x_0,\ldots,x_{n-1})$

To accomplish this, we will simulated a Markov chain with state space
$\X^n$ whose equilibrium distribution is
$\pi(x_0,\ldots,x_{n-1})$. The state at iteration $i$ of this chain
will be written as $x^{(i)} = (x_0^{(i)},\ldots,x_{n-1}^{(i)})$.  The
transition probabilities for this Markov chain will be denoted using
$Q$.  In particular, we will use some initial distribution for the
state, $Q(x^{(0)})$, and will simulate the chain according to the
transition probabilities $Q(x^{(i)}\,|\,x^{(i-1)})$.  For validity of
the sampling method, we need these transitions to leave $\pi$
invariant: 
\beq
  \pi(x') & = & \sum_{x} \pi(x) Q(x'\,|\,x),\ \ \ \mbox{for all $x'$ in $\X^n$}
\eeq
(If $\X$ is continuous, the sum is replaced by an integral.)  This 
is implied by the detailed balance condition:
\beq
 \pi(x) Q(x'\,|\,x) & = & \pi(x') Q(x\,|\,x'),\ \ \ 
    \mbox{for all $x$ and $x'$ in $\X^n$}
\eeq

The transition $Q(x^{(i)}\,|\,x^{(i-1)})$ is defined using a set of auxiliary
Markov chains, one for each time step, whose state spaces are $\X$, and 
whose transition probabilities, written as $R_t(\cdot\,|\,\cdot)$, leave
a specified ``pool'' distribution, $\rho_t$, invariant.
The transitions for the reversal of this chain with respect to 
$\rho_t$ will be denoted by $\tilde R(\cdot\,|\,\cdot)$.
These transitions satisfy the following condition:
\beq
 \rho_t(x) R_t(x'\,|\,x) & = & \rho_t(x') \tilde R_t(x\,|\,x'),\ \ \ 
  \mbox{for all $x$ and $x'$ in $\X$}
\eeq
Note that if the transitions $R_t$ satisfy detailed balance with respect to 
$\rho_t$, $\tilde R_t$ will be the same as $R_t$.

For each time, $t$, the transitions $R_t$ and $\tilde R_t$ are used 
to produce a pool of
$K$ candidate states, $\C_t$, one of which is the current state, $x^{(i-1)}_t$.
The new sequence, $x^{(i)}$, is randomly selected from 
among all sequences whose states at each time $t$ are in $\C_t$, using 
a form of the forward-backward procedure.

In detail, the pool of candidate states for time $t$ is found as 
follows:\vspace*{-4pt}
\begin{enumerate}
  \item[1)] Pick an integer $J_t$ uniformly from $\{0,\ldots,K-1\}$.

  \item[2)] Let $x_t^{[0]} = x_t^{(i-1)}$. 

  \item[3)] For $j$ from $1$ to $J_t$, randomly pick $x_t^{[j]}$ 
             according to
             the transition probabilities $R_t(x_t^{[j]}\,|\,x_t^{[j-1]})$.

  \item[4)] For $j$ from $-1$ down to $-K+J_t+1$, randomly pick $x_t^{[j]}$ 
             according to the reversed
             transition probabilities, $\tilde R_t(x_t^{[j]}\,|\,x_t^{[j+1]})$.
 
  \item[5)] Let $\C_t$ be the pool consisting of $x^{[j]}_t$, for $j \in 
            \{ -K+J_t+1, \ldots, 0, \ldots, J_t\}$.  If some of the $x^{[j]}_t$
            are the same, they will be present in the pool more than once.
\end{enumerate}

Once the pools of candidate states have been found, a new state sequence, 
$x^{(i)}$, is picked from among all sequences, $x$,
for which every $x_t$ is in $\C_t$.  The probability of picking $x$ is
proportional to $\pi(x)/\prod\rho_t(x_t)$, which is proportional to
\beq
{ P(x_0) \prod_{t=1}^{n-1} P(x_t\,|\,x_{t-1}) \prod_{t=0}^{n-1} P(y_t\,|\,x_t)
   \over \prod_{t=0}^{n-1}\rho_t(x_t) }\label{eq-poolpr}
\eeq
If duplicate states occur in some of the pools, they are 
treated as if they were distinct when picking a sequence in this way.
In effect, we pick indexes of states in these pools, with probabilities
as above, rather than states themselves.  The distribution of these
sequences of indexes can be regarded as the posterior distribution for
a hidden Markov model, with the transition probability
from state $j$ at time $t-1$ to state $k$ at time $t$ being
proportional to $P(x^{[k]}_t\,|\,x^{[j]}_{t-1})$, and the 
probabilities of the hypothetical observed symbols being proportional
to the remaining factors above, $P(y_t\,|\,x^{[k]}_t) / \rho_t(x^{[k]}_t)$.
Crucially, it is possible, using the forward-backward technique, to
randomly pick a new state from this distribution in time growing
linearly with $n$, even though the number of possible sequences grows
as $K^n$.

\section{\hspace*{-7pt}Proof of Correctness}\label{sec-proof}\vspace*{-10pt}

To show that a Markov chain with these transitions will converge to
$\pi$, we need to show that it leaves $\pi$ invariant, and that the
chain is ergodic.  Ergodicity need not always hold, and proving that
it does hold may require considering the particulars of the model.
However, it is easy to see that the chain will be ergodic if all
possible state sequences have non-zero probability density under
$\pi$, the pool distributions, $\rho_t$, have non-zero density
everywhere, and the transitions $R_t$ are ergodic.  This probably
covers most problems that arise in practice.

To show that the transitions $Q(\cdot\,|\,\cdot)$ leave $\pi$ invariant,
it suffices to show that they satisfy detailed balance with respect to
$\pi$.  This will follow from the stronger condition that the probability
of moving from $x$ to $x'$ (starting from a state picked from $\pi$) 
with given values for the $J_t$ and given pools of candidate states, $\C_t$, 
is the same as the corresponding probability of moving from $x'$ to $x$
with the same pools of candidate states and with values $J'_t$ defined
by $J'_t = J_t - h_t$, where $h_t$ is the index (from $-K+J_t+1$ to $J_t$)
of $x'_t$ in the candidate pool.

The probability of such a move from $x$ to $x'$ is the product of
several factors.  First, there is the probability of starting from $x$
under $\pi$, which is $\pi(x)$.  Then, for each time $t$, there is the
probability of picking $J_t$, which is $1/K$, and of then producing the
states in the candidate pool using the transitions $R_t$ and $\tilde R_t$,
which is
\beq
   \lefteqn{\prod_{j=1}^{J_t} R_t(x_t^{[j]}\,|\,x_t^{[j-1]})\ \ \times\!\!\!
   \prod_{j=-K+J_t+1}^{-1} \tilde R_t(x_t^{[j]}\,|\,x_t^{[j+1]})}\ \ \ \ \ \ \
   \nonumber\\[5pt]
   & = & 
   \prod_{j=0}^{J_t-1} R_t(x_t^{[j+1]}\,|\,x_t^{[j]})\ \ \times\!\!\!
   \prod_{j=-K+J_t+1}^{-1}\!\! R_t(x_t^{[j+1]}\,|\,x_t^{[j]})\
                    {\rho_t(x_t^{[j]}) \over \rho_t(x_t^{[j+1]})} \\[5pt]
   & = & 
   {\rho_t(x_t^{[-K+J_t+1]}) \over \rho_t(x_t^{[0]})}\!
   \prod_{j=-K+J_t+1}^{J_t-1}\!\! R_t(x_t^{[j+1]}\,|\,x_t^{[j]}) 
\eeq
Finally, there is the probability of picking $x'$ from among the sequences
with states from the pools, $\C_t$, which is proportional to 
$\pi(x')/\prod\rho_t(x'_t)$.  The product of all these factors is
\beq
\lefteqn{\pi(x) \ \times\  (1/K)^n \ \times\  \prod_{t=0}^{n-1} \left[ 
   {\rho_t(x_t^{[-K+J_t+1]}) \over \rho_t(x_t^{[0]})}\!
   \prod_{j=-K+J_t+1}^{J_t-1}\!\! R_t(x_t^{[j+1]}\,|\,x_t^{[j]}) \right]
   \ \times\ {\pi(x') \over \prod_{t=0}^{n-1}\rho_t(x'_t)}}
   \ \ \ \ \ \ \ \ \nonumber\\[5pt]
   & = &
   (1/K)^n\ {\pi(x)\pi(x')\over \prod_{t=0}^{n-1} \rho(x_t)\rho(x'_t)}
   \ \prod_{t=0}^{n-1} \left[ \rho_t(x_t^{[-K+J_t+1]})\!\!
   \prod_{j=-K+J_t+1}^{J_t-1}\!\! R_t(x_t^{[j+1]}\,|\,x_t^{[j]})\right] \ \ \ \
\eeq
The corresponding expression for a move from $x'$ to $x$ is identical,
apart from a relabelling of candidate state $x_t^{[j]}$ as $x_t^{[j-h_t]}$.

\section{\hspace*{-7pt}An Example Class of Models}\label{sec-ex}\vspace*{-10pt}

As a simple concrete example, consider a model in which the state space $\X$
and the observation space, $\Y$, are both $\Re$.  Let each observation
be simply the state plus Gaussian noise of standard deviation $\sigma$ --- ie,
$P(y_t\,|\,x_t) = N(y_t\,|\,x_t,\,\sigma^2)$ --- and let the state
transitions be defined by $P(x_t\,|\,x_{t-1}) = N(x_t\,|\,\tanh(\eta x_{t-1}),\,
\tau^2)$, for some constant expansion factor $\eta$ and transition noise
standard deviation $\tau$.  

Let us choose the pool distributions, $\rho_t$, to be normal, with
some means $\mu_t$ and standard deviations $\nu_t$, which may depend
on $y_0,\ldots,y_{n-1}$, but not on $x_0,\ldots,x_{n-1}$.  For
example, we might fix $\mu_t=0$ and $\nu_t=1$ for all $t$, or we might
let $\rho_t$ be the posterior distribution for $x_t$ given $y_t$,
based on an improper flat prior, so that $\mu_t=y_t$ and
$\nu_t=\sigma$, or we might let $\rho_t$ be some more elaborate
approximation to the marginal distribution of $x_t$ given
$y_0,\ldots,y_{n-1}$.

Of the many transitions that would leave $\rho_t$ invariant, we might
choose $R_t$ to be of the following form:
\beq
  R_t(x'\,|\,x) & = & 
    N(x'\,|\,\mu_t+\alpha (x\!-\!\mu_t),\,(1\!-\!\alpha^2)\nu_t^2)
\label{eq-Rex}\eeq
where $\alpha$ is an adjustable parameter in $(-1,+1)$.  When $\alpha=0$,
the states in the pool (other than the current state) are drawn independently 
from $\rho_t$.  These transitions satisfy detailed balance with respect to 
$\rho_t$, so $\tilde R_t$ is the same as $R_t$.

The forward-backward algorithm will pick a state
sequence from among those that can be constructed using states from 
the candidate pools, with probabilities given by equation~(\ref{eq-poolpr}).  
In the particular case when $\mu_t=y_t$ and $\nu_t^2=\sigma^2$, for
which $\rho_t(x_t)$ is proportional to $P(y_t\,|\,x_t)$, these
probabilities simplify to being proportional to
\beq
 P(x_0) \prod_{t=1}^{n-1} P(x_t\,|\,x_{t-1})
\eeq
Note that
despite appearances, this distribution cannot be sampled from using
a forward pass alone, since $P(x_t\,|\,x_{t-1})$ need not sum to one
for $x_t$ in $\C_t$.

\section{\hspace*{-7pt}Demonstration}\label{sec-demo}\vspace*{-10pt}

The characteristics of the state and observation sequences produced
using the models of the previous section vary considerably with the
choice of $\sigma$, $\eta$, and $\tau$.  For some choices, simple
forms of the Metropolis algorithm that update each $x_t$ separately
can perform better than the embedded HMM method, since these simple
methods have lower overhead.  Here I will demonstrate that the
embedded HMM can perform better than such single-state updating
methods when the states are highly dependent.

\begin{figure}[p]

\vspace*{-33pt}

\centerline{\psfig{file=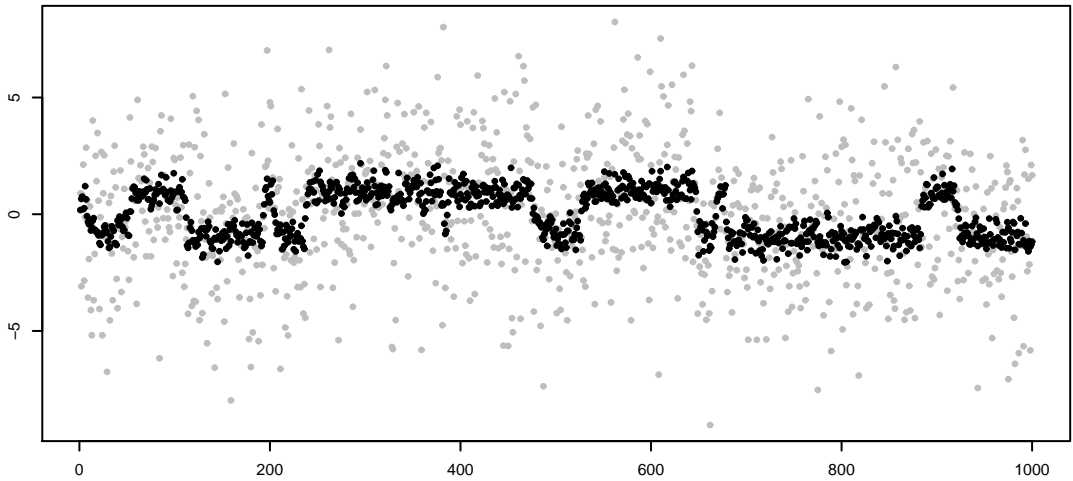}}

\vspace*{-27pt}

\caption[]{A state sequence (black dots) and observation sequence (gray dots)
of length 1000 produced by the model with $\sigma=2.5$, $\eta=2.5$, and 
$\tau=0.4$.}\label{fig-seqb}

\end{figure}

Figure~\ref{fig-seqb} shows a sequence, $x_0,\ldots,x_{n-1}$, and
observation sequence, $y_0,\ldots,y_{n-1}$, produced using
$\sigma=2.5$, $\eta=2.5$, and $\tau=0.4$, with $n=1000$.  The state
sequence stays in the vicinity of $+1$ or $-1$ for long periods, with
rare switches between these regions.  Because of the large observation
noise, there is considerable uncertainty regarding the state sequence
given the observation sequence, with the posterior distribution
assigning fairly high probability to sequences that contain short-term
region switches that are not present in the actual state sequence, or
that lack some of the short-term switches that are actually present.
It is difficult for a method that updates only one state at a time to
explore such a posterior distribution, because it must move through
low-probability intermediate states in which a switch to the opposite
region is followed immediately by a switch back.

\begin{figure}[p]

\vspace*{-22pt}

\centerline{\psfig{file=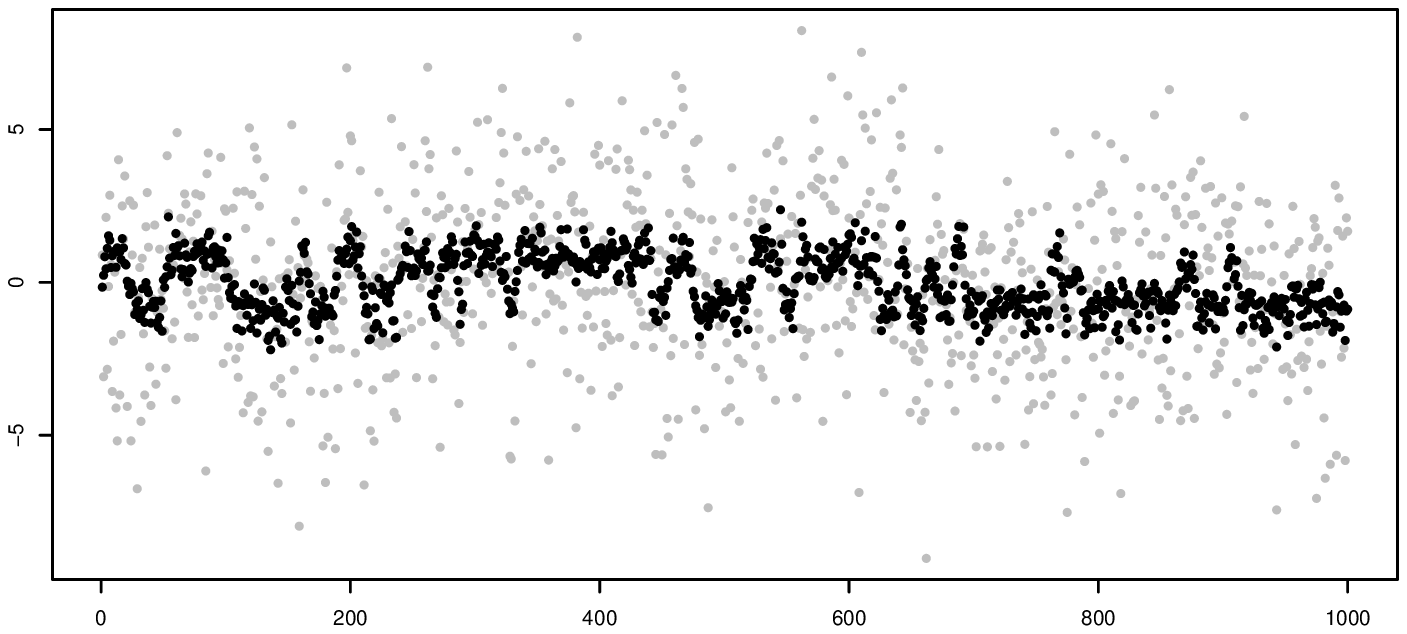}}

\vspace*{-40pt}

\centerline{\psfig{file=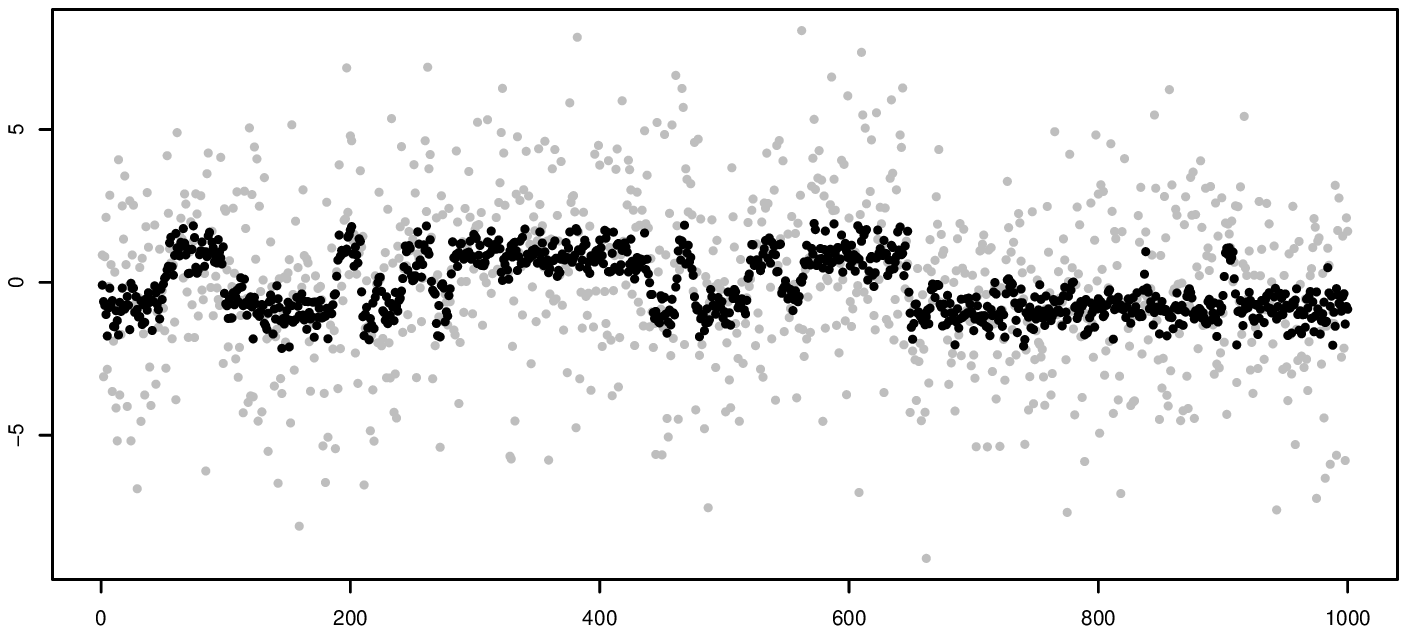}}

\vspace*{-27pt}

\caption[]{State sequences (black dots) produced after one embedded HMM 
update (top) and two updates (bottom), starting with the states set 
equal to the data points (gray dots), for the same model and data
as Figure~\ref{fig-seqb}.  The embedded HMM used $K=10$, $\mu_t=0$,
$\nu_t=1$, and $\alpha=0$.}\label{fig-hmmb12}

\end{figure}

Figure~\ref{fig-hmmb12} shows that embedded HMM sampling works well
for this problem, using $K=10$ states and the simple choice of
$\mu_t=0$ and $\nu_t=1$ for the pool distributions, and $R_t$ as in
equation~(\ref{eq-Rex}), with $\alpha=0$.  We can see that only two
updates produce a state sequence with roughly the correct
characteristics.  

\begin{figure}[p]

\vspace*{-33pt}

\centerline{\psfig{file=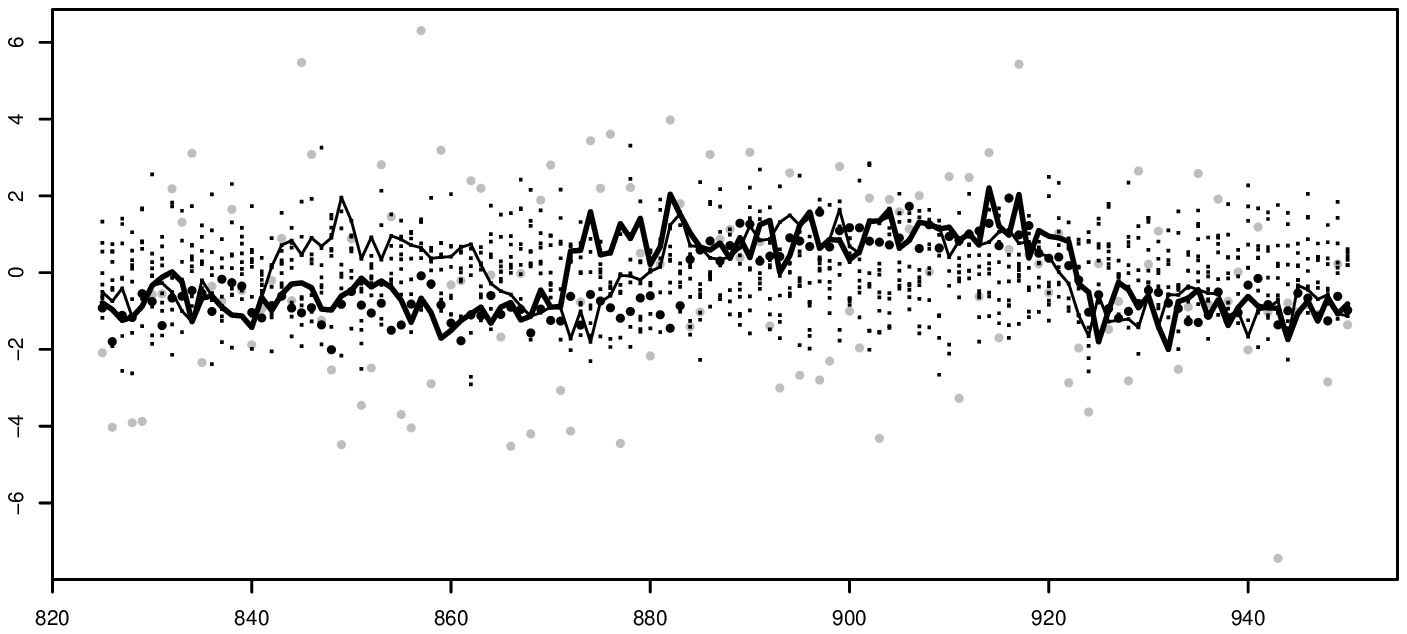}}

\vspace*{-27pt}

\caption[]{Closeup of an embedded HMM update.  The true state sequence 
is shown by black dots and the observation sequence by gray dots.  The
current state sequence is shown by the dark line.  The pools of 
states used for the update are shown as small dots, and the new state
sequence picked by the embedded HMM by the light line.
}\label{fig-hmmbp}

\end{figure}

Figure~\ref{fig-hmmbp} demonstrates how a single embedded HMM update
can make a large change to the state sequence.  It shows a portion of
the state sequence after 99 updates, the pools of states produced for
the next update, and the state sequence found by the embedded HMM
using these pools.  A large change is made to the state sequence in
the region from time 840 to 870, with states in this region switching
from the vicinity of $-1$ to the vicinity of $+1$.

\begin{figure}[p]

\vspace*{-33pt}

\centerline{\psfig{file=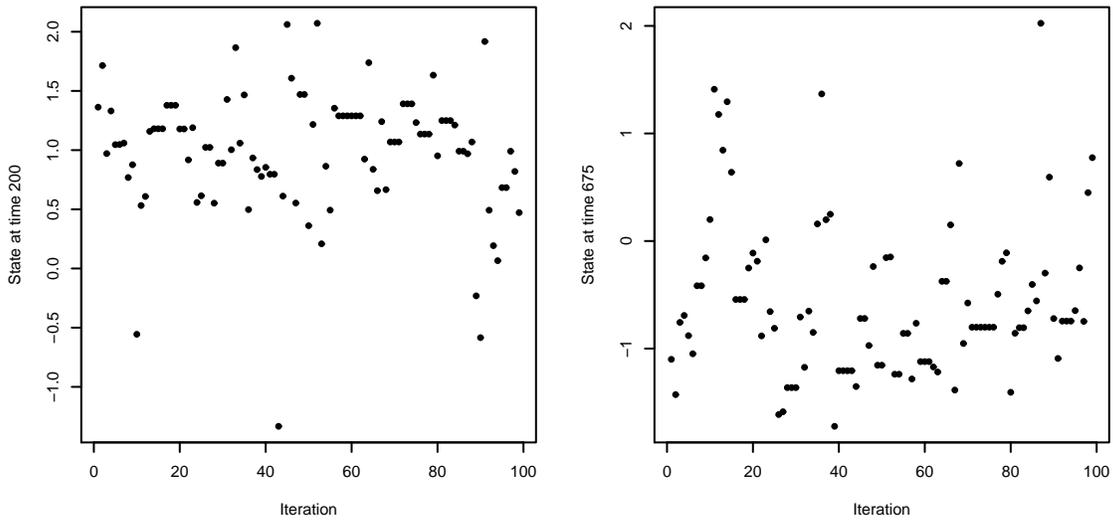}}

\vspace*{-15pt}

\caption[]{Traces of states during an embedded HMM run.  The left plot 
shows the state at time 200 after each of the first 99 updates; the
right plot shows the same for the state at time 675.}\label{fig-hmmbt}

\end{figure}

In Figure~\ref{fig-hmmbt}, the state at two time points is plotted
over the course of 99 embedded HMM updates.  Both points correspond to
short-term switches in the actual state sequence.  In the posterior
distribution, there is uncertainty about the true state at these
points, with non-negligible probability for values near $-1$ and for
values near $+1$.  We see in both plots that the embedded HMM moves
between these two regions.

\begin{figure}[p]

\vspace*{-22pt}

\centerline{\psfig{file=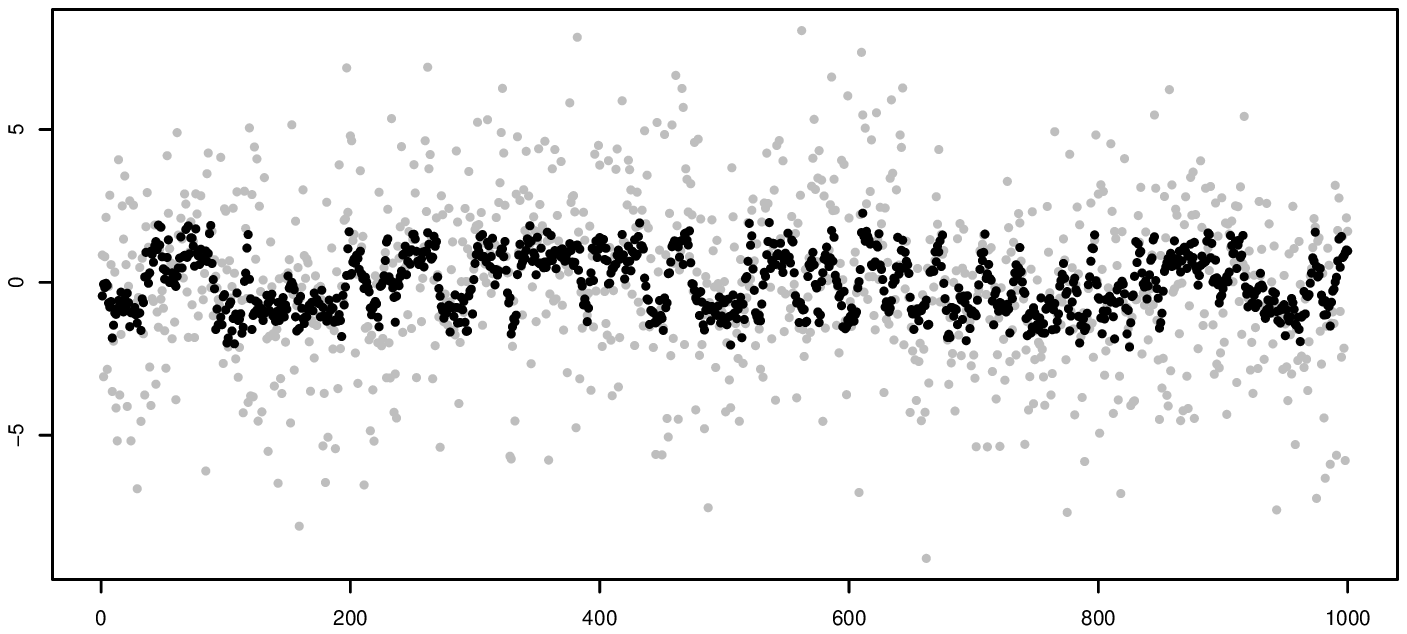}}

\vspace*{-40pt}

\centerline{\psfig{file=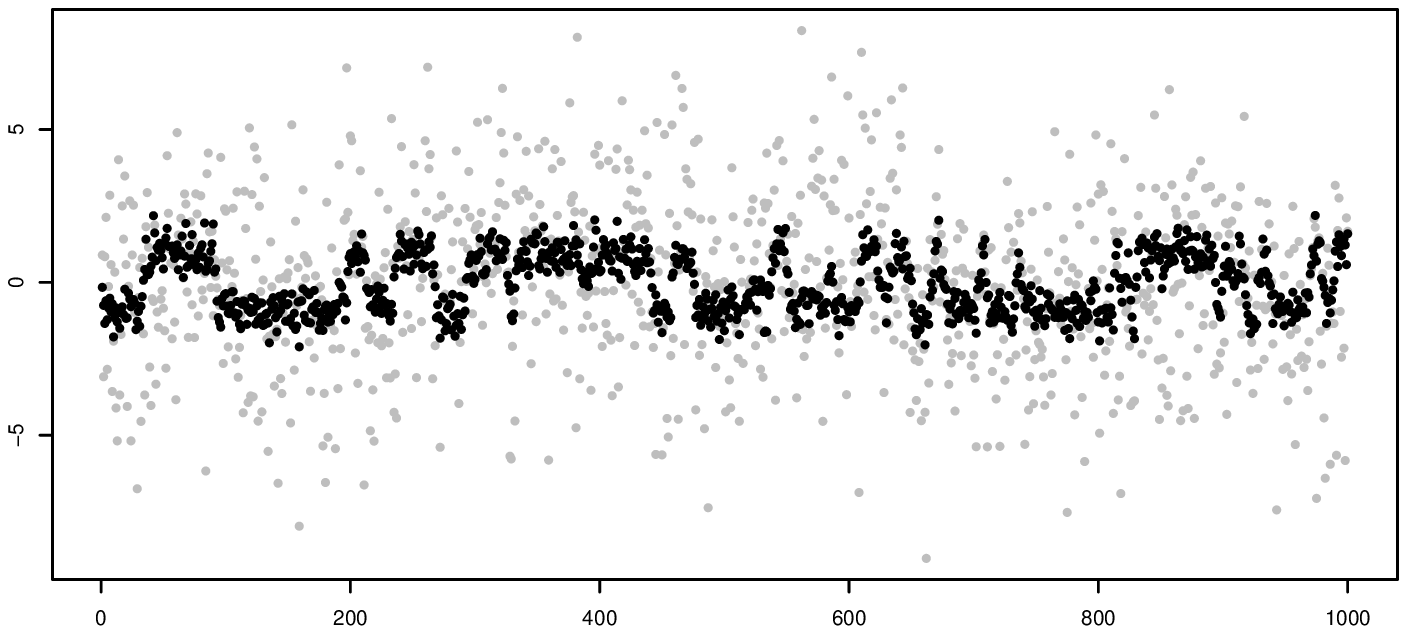}}

\vspace*{-27pt}

\caption[]{State sequences (black) produced after 50 single-state
Metropolis updates (top) and after 100 updates (bottom), starting with the 
states set equal to the data points (gray).}\label{fig-metb50100}

\end{figure}

\begin{figure}[p]

\vspace*{-23pt}

\centerline{\psfig{file=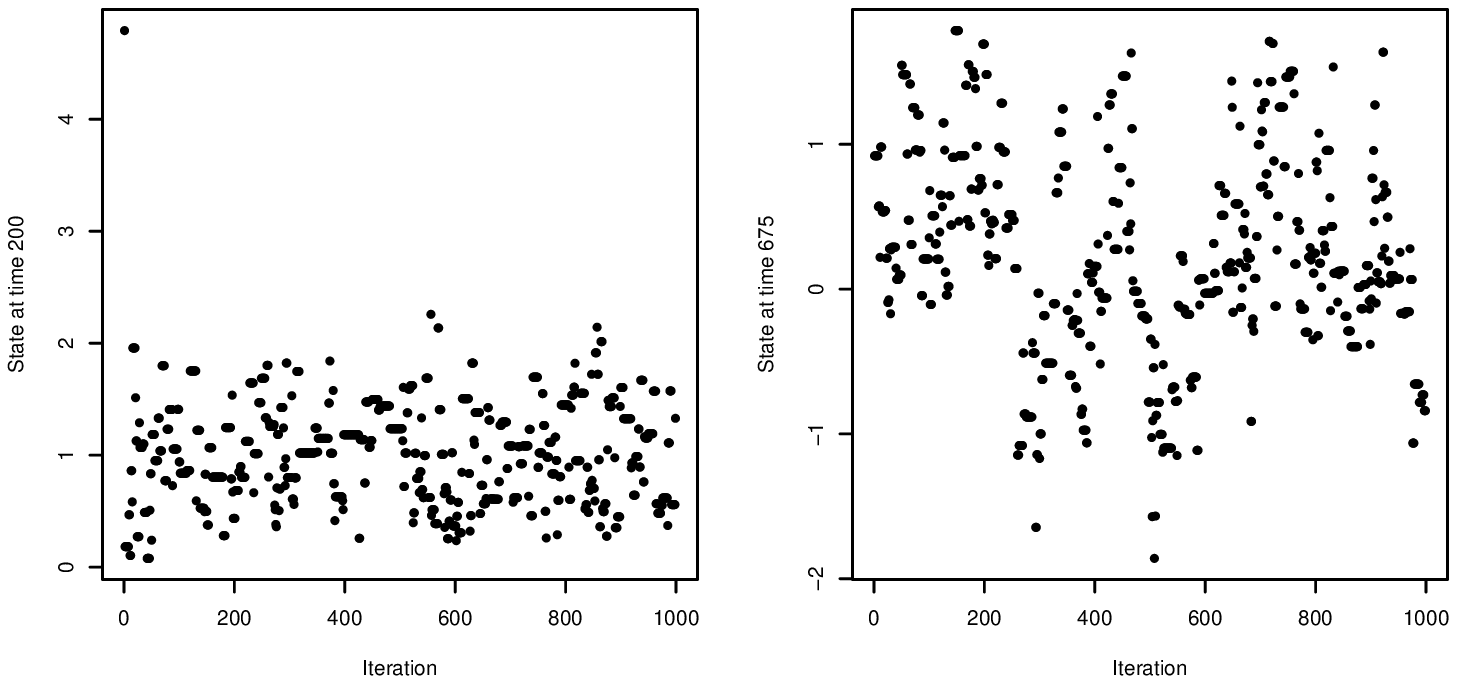}}

\vspace*{-15pt}

\caption[]{Traces of states during a Metropolis run.  The left plot 
shows the state at time 200 after each of the first 999 updates; the
right plot shows the same for the state at time 675.}\label{fig-metbt}

\end{figure}

In contrast, simple Metropolis methods that update one state at a time
do much less well for this problem.  Figure~\ref{fig-metb50100} shows
the state sequences produced after 50 and 100 iterations of a
Metropolis method in which each iteration updates each state in turn,
using a $N(0,1)$ proposal distribution.  Even after 100 iterations,
the state sequence does not closely resemble the actual state sequence
(in particular, it contains too many short-term switches).  The traces
of states at times 200 and 675 in Figure~\ref{fig-metbt} confirm that
these Metropolis updates do not move around the posterior distribution
efficiently.  The state at time 200 never reaches the vicinity of $-1$
during these 999 iterations.  The state at time 675 does visit the
vicinity of both $-1$ and $+1$, but the values show very high
autocorrelations.  Simple Metropolis updates with various other
proposal distributions performed similarly, or worse.

On the other hand, one iteration of these simple Metropolis methods is
approximately 30 times faster than one of the embedded HMM updates
(with $K=10$), when both methods are implemented in the interpretive R
language.  In this example, however, the greater efficiency of the
embedded HMM updates more than outweighs this.

With other settings of the $\sigma$, $\eta$, and $\tau$ parameters,
different pool distributions are preferable to the simple $N(0,1)$
distribution used for this demonstration.  In particular, letting the
pool distribution for $x_t$ depend on $y_t$ or on a window of
observations in its vicinity is sometimes better.  I have not found
setting $\alpha$ to a non-zero value to be beneficial for this model,
but I expect that setting $\alpha$ close to one in order to produce a
pool of states in the vicinity of the current state will be useful in
higher-dimensional problems.

\section*{Acknowledgement}\vspace{-10pt}

I thank Sam Roweis and Matthew Beal for helpful discussions.  This
research was supported by the Natural Sciences and Engineering
Research Council of Canada.

\section*{Reference}\vspace*{-10pt}

\leftmargini 0.2in
\labelsep 0in

\begin{description}
\itemsep 0pt

\item

  Scott, S.\ L.\ (2002) ``Bayesian methods for hidden Markov
  models: Recursive computing in the 21st century'', \textit{Journal 
  of the American Statistical Association}, vol.~97, pp.~337--351.

\end{description}

\end{document}